\documentclass[a4paper,11 pt]{amsart}

\RequirePackage{amsmath, amssymb, amsthm, amsfonts}
\usepackage[all]{xy}
\usepackage{enumitem}
\usepackage{hyperref}
\usepackage{fullpage}
\usepackage{hyperref}
\usepackage[mathscr]{eucal}
\usepackage[active]{srcltx}
\newtheorem{theo}{Theorem}[section]

\newtheorem{prop}[theo]{Proposition}

\theoremstyle{definition}

\newtheorem{quest}[theo]{Question}
\newtheorem{rem}[theo]{Remark}

\newcommand{\oo}{\mathcal{O}}

\newcommand{\EE}{\mathscr{E}}
\newcommand{\mL}{\mathcal{L}}
\newcommand{\mQ}{\mathcal{Q}}

\newcommand{\SL}{{\mathrm{SL}}}

\newcommand{\GL}{{\mathrm{GL}}}

\newcommand{\Alb}{\textrm{Alb}}

\newcommand{\ZZ}{\mathbb{Z}}

\newcommand{\wA}{\widehat{A}}
\newcommand{\wS}{\widehat{S}}

\newcommand{\lr}{\longrightarrow}
\newcommand{\LL}{\mathscr{L}}

\title{On the classification of surfaces of general type with $p_g=q=2$}
\author{Matteo Penegini}
\address{Matteo Penegini\\
Dipartimento di Matematica ``F. Enriques" \\ Via Saldini 50, I-20100
Milano, Italy} \email{matteo.penegini@unimi.it}

\begin{document}


\maketitle

\tableofcontents





\section{Introduction}\label{sec.intro}

The following is an extended version of the talk which I gave at
the ``\textit{XIX Congresso dell'UMI}'' in Bologna in September
2011. The aim of this paper is twofold: first, to give an overview on the recent development in
the classification of surfaces of general type with $p_g=q=2$; second, to point out some of the problems that are still open.
Some of the results that will appear in this article have been
proven in collaboration with Francesco Polizzi.

We shall use the standard notation from the theory of complex
algebraic surfaces. Let $S$ be a smooth, complex, projective
\emph{minimal} surface $S$ of \emph{general type}; this means that
the canonical divisor $K_S$ of $S$ is \emph{big} and
\emph{nef}.



The principal numerical invariants for the study of minimal
surfaces of general type are

\begin{itemize}
    \item the \emph{geometric genus} $p_g(S):=h^0(S, \,
\Omega^2_S)=h^0(S, \, \mathcal{O}_S(K_S))$,
    \item the
\emph{irregularity} $q(S):=h^0(S, \, \Omega^1_S)$, and
    \item the \emph{self
intersection of the canonical divisor} $K^2_S$.
\end{itemize}

As a matter of fact, these determine all the other classical
invariants, as
\begin{itemize}
   \item the
\emph{Euler-Poincar\'e characteristic}
$\chi(\mathcal{O}_S)=1-q(S)+p_g(S)$,
    \item the \emph{topological Euler number}
    $e(S)=12\chi(S)-K^2_S$, and
    \item the \emph{plurigenera} $P_n(S)=\chi(S)+(^n_2)K^2_S$.
\end{itemize}

By a theorem of Bombieri, a minimal surface of general type $S$
with fixed invariants is birationally mapped to a normal surface
$X$ in a fixed projective space of dimension $P_5(S)-1$. Moreover,
$X$ is uniquely determined and is called the \emph{canonical
model} of $S$. Let us recall Gieseker's Theorem (see \cite{Gie77}).

\begin{theo}
 There exists a quasi-projective coarse moduli space
$\mathcal{M}_{K^2_S, \chi}$ for canonical models of surfaces of
general type $S$ with fixed invariants $K^2_S$ and $\chi$.
\end{theo}

In particular, we can also consider the subscheme
$\mathcal{M}_{K^2_S,p_g,q} \subset \mathcal{M}_{K^2_S,\chi}$ of surfaces with fixed $K^2_S$, $p_g(S)=p_g$ and $q(S)=q$, which is a union of connected components of $\mathcal{M}_{K^2_S,\chi}$. The
first basic question we would like to answer is

\begin{quest}{\it For
which values of $(K^2_S,p_g,q)$ is $\mathcal{M}_{K^2_S,p_g,q}$ not
empty}?
\end{quest}
Much is known about this question; indeed, we have the following classical
inequalities:

\begin{itemize}
\item $K^2_S \geq 1$ and $\chi \geq 1$.
\medskip
\item $K^2_S \leq 9 \chi(S)$ (Bogomolov-Miyaoka-Yau).
\medskip
\item $K^2_S \geq 2 p_g$, if $q>0$, (Debarre).
\medskip
\item $K^2_S \geq 2\chi - 6$ (Noether).
\end{itemize}
Once we have established that $\mathcal{M}_{K^2_S,p_g,q} \neq
\emptyset$, the second basic question is
\begin{quest}\textit{Can we describe
$\mathcal{M}_{K^2_S,p_g,q}$?}
\end{quest}
For example, can we find out the number of
connected components it consists and their dimensions?

We are interested in surfaces with small invariants, and in
particular in surfaces with $\chi(S)=1$, which is the smallest
possible value for a surface of general type. This implies that
$p_g=q$ and the above inequalities yield
\[ 0 \leq p_g \leq 4 .
\]

A partial classification of these surfaces has already
been accomplished. In the case $\mathbf{p_g=q=4}$ we have a full
classification theorem.

\begin{theo}\cite{Be82} If $S$ is a minimal surface of general type with $p_g=q=4$, then $S$ is a product of two curves of genus
$2$ and $K^2_S=8$. Moreover, $\mathcal{M}_{8,4,4}$ consists of
exactly one connected component of dimension $6$.
\end{theo}

Also, in the case $\mathbf{p_g=q=3}$ we have a full classification
theorem due to several authors: \cite{CCML}, \cite{HP02}, and
\cite{Pi02}.

\begin{theo}\label{theo_pgq3} Let $S$ be a
minimal surface of general type with $p_g=q=3$, then there are
only two possibilities:
\begin{enumerate}
\item $K^2_S=8$, and $S = (C \times F)/G$ where $C$ is a curve of
genus $2$, $F$ is a curve of genus $3$ and $G \cong
\mathbb{Z}/2\ZZ$. Here $G$ acts freely and diagonally on the product $C \times F$, on $C$ as a hyperelliptic involution, and on $F$ as a fixed
point free involution. Moreover, $\mathcal{M}_{8,3,3}$ consists of
exactly one connected component of dimension $5$.

\item $K^2_S=6$, and $S$ is the symmetric square of a genus $3$
curve.  Moreover, $\mathcal{M}_{6,3,3}$ consists of exactly one
connected component of dimension $6$.
\end{enumerate}
\end{theo}

On the other hand, in the case $p_g \leq 2$ a complete classification
theorem is still missing. It seems that the classification
becomes more complicated as the value of $p_g$ decreases. In this
paper, we address the case $p_g=q=2$.

Since we are dealing with irregular surfaces, i.e., with $q > 0$,
a useful tool that we can use is the Albanese map. The
\emph{Albanese variety} of $S$ is defined as
$\Alb(S):=H^0(\Omega^1_S)^{\vee}/H_1(S,\ZZ)$. By Hodge theory,
$\Alb(S)$ is an abelian variety. For a fixed base point $x_0 \in
S$, we define the \emph{Albanese} morphism
\[ \alpha_{x_0}\colon S \longrightarrow \Alb(S), \quad \quad x \mapsto
\int^x_{x_0}.
\]
If we choose a different base point in $S$, the Albanese morphism
changes by a translation of $\Alb(S)$; so we often ignore the base
point and write $\alpha$.

The dimension of $\alpha(S)$ is called the \emph{Albanese
dimension} of $S$ and it is denoted by $\textrm{Albdim(S)}$. If
$\textrm{Albdim(S)}=2$, we say that $S$ has \emph{maximal Albanese
dimension}. If, moreover, $q(S)>2$, we say that $S$ is of
\emph{Albanese general type}.

For surfaces with $q(S)=2$, we have two possibilities:
\begin{enumerate}
\item either $\alpha(S)$ is a curve, or
\item  $\alpha$ is a generically finite
cover of an abelian surface.
\end{enumerate}

In the latter case, the degree of $\alpha$ plays a crucial role.
Indeed, if $q(S)=2$, then the degree of the Albanese map $\alpha$
is a topological invariant, see e.g. \cite[Section 5]{Ca11}.

In this paper, we want to present the up-to-date list of the known
surfaces of general type with $\mathbf{p_g=q=2}$, which can be
shortly summarized in the following table.
\begin{center}
\begin{tabular}{|c|c|c|c|c|c|c|}
  \hline
n. &$K^2_S$ & $Albdim(S)$ & $\textrm{deg}{\alpha}$ &  $\sharp$ Families & $dim$ & Name \\
  \hline
1 & $8$ & $1$ & $-$ &  $24$ & $3^{15},4^6,5^2,6$ & Isog. to a Prod. \\ \hline
2 & $8$ & $2$ & $ \leq 6$ &  $4$ &  $3^3,4$     & Isog. to a Prod. \\ \hline
3 & $6$ & $2$ & $4$ &  $1$ &  $4$ &  \\ \hline
4 & $6$ & $2$ & $2$ &  $3$ &  $4^2,3$ &  \\ \hline
5 & $5$ & $2$ & $3$ &  $1$ &  $4$ & Chen-Hacon surf. \\ \hline
6 & $4$ & $2$ & $2$ &  $1$  &  $4$   & 
\\ \hline
\end{tabular}
\end{center}
\begin{center}
Table 0.
\end{center}
In the table we used the notation $3^{15}$ to mean that
there are $15$ families of dimensions $3$.

The paper is organized as follows. We shall describe the families
of our surfaces following the numeration of Table 0.

Indeed, in the first section we shall treat the case where
$\alpha(S)$ is a curve. Here we have a complete classification
theorem. We introduce the surfaces isogenous to a product of
curves and some of their properties. The end of the section is
dedicated to product-quotient surfaces, which are a generalization
of surfaces isogenous to a product.

In the second section, we describe the four families of surfaces
with $p_g=q=2$ and $K^2_S=6$. Three families are characterized by
the fact that their elements have Albanese map of degree $2$,
while the last one has Albanese map of degree $4$.

The third section is devoted to describe \emph{Chen-Hacon
surfaces}. These surfaces were introduced by Chen and Hacon in
\cite{CH06} and together with F. Polizzi were able to give a description of the
connected component of the moduli space of surfaces of general type
they belong to.

The last section is devoted to the last cases and
to give an account of the open problems on this topic.

\bigskip

 \textbf{Acknowledgments.} The author is grateful to G. Bini, F. Catanese, A. Garbagnati, E.C. Mistretta, R. Pignatelli e F. Polizzi for reading and commenting the paper. 

%

\section{Product-Quotient Surfaces with $p_g=q=2$}

%

A surface $S$ is said to be \emph{isogenous to a (higher) product of
curves} if and only if $S$ is a quotient $S:=(C \times F)/G$, where
$C$ and $F$ are curves of genus at least two, and $G$ is a finite
group acting freely on $C \times F$.

Let $S$ be a surface isogenous to a higher product, and
$G^{\circ}:=G \cap(Aut(C) \times Aut(F))$. Then $G^{\circ}$ acts
on the two factors $C$ and $F$ and diagonally on the product $C
\times F$. If $G^{\circ}$ acts faithfully on both curves, we say
that $S= (C \times F)/G$ is a \emph{minimal realization} of  $S$.
In \cite{cat00}, the author proves that any surface isogenous to a higher
product admits a unique minimal realization. From now on we shall
work only with minimal realizations.

There are two cases: the \emph{mixed} case where the action of $G$
exchanges the two factors (in this case $C$ and $F$ are isomorphic
and $G^{\circ} \neq G$); the \emph{unmixed} case (where
$G=G^{\circ}$, and therefore it acts diagonally).

Moreover, we observe that a surface isogenous to a product
of curves is of general type. It is always minimal and its
numerical invariants are explicitly given in terms of the genera
of the curves and the order of the group. Indeed, we have the
following proposition.
\begin{prop} Let $S=(C \times F)/G$ be a surface isogenous to a higher product of curves, then:
\begin{equation*}\label{eq.chi.isot.fib}
\chi(S)=\frac{(g(C)-1)(g(F)-1)}{|G|},
\end{equation*}
\begin{equation*}\label{eq.euler.isot.fib}
e(S)=\frac{4(g(C)-1)(g(F)-1)}{|G|},
\end{equation*}
\begin{equation*}\label{eq.k2.isot.fib}
K^2_S=\frac{8(g(C)-1)(g(F)-1)}{|G|}.
\end{equation*}

\end{prop}

In the unmixed case  $G$ acts separately on $C$ and $F$, and the
two projections $\pi_C \colon C \times F \lr C, \; \pi_F \colon C
\times F \lr F$ induce two isotrivial fibrations $\alpha \colon S
\lr C/G, \; \beta \colon S \lr F/G$, whose smooth fibres are
isomorphic to $F$ and $C$, respectively.

A surface isogenous to a product of unmixed type $S:=(C \times
F)/G$ is said to be of \emph{generalized hyperelliptic type} if:
\begin{enumerate}
\item the Galois covering $\pi_1: C \rightarrow C/G
    $ is unramified, and
\item the quotient curve $F/G$ is isomorphic to
    $\mathbb{P}^1$.
\end{enumerate}

Surfaces of generalized hyperelliptic type play a crucial role
among surfaces with $p_g=q=2$ because of the following
proposition.

\begin{prop}[\cite{zucc} Prop. 4.2]\label{Zucco} If $S$ is a surface of general type
with $p_g=q=2$ and $\textrm{Albdim(S)}=1$. Then $S$ is of
generalized hyperelliptic type.
\end{prop}

By this proposition, once we classify all the surfaces
isogenous to a product with $p_g=q=2$, we also achieve the
classification of all the ones such that $\alpha(S)$ is a curve.
This work was started by Zucconi in \cite{zucc} and
completed in \cite{Pe11}. We can summarize the results in the
following theorem.

\begin{theo}[\cite{Pe11} Theorem 1.1]\label{zerozero}
Let $S$ be a surface isogenous to a higher product of curves with
$p_g=q=2$. Then we have the following possibilities:
\begin{enumerate}
\item If $\textrm{Albdim}(S))=1$, then $S \cong (C \times F)/G$
and it is of generalized hyperelliptic type. The classification of
these surfaces is given by the cases labelled with GH in Table 1,
where we specify the possibilities for the genera of the two
curves $C$ and $F$, and for the group $G$.

\item If $\textrm{Albdim}(S))=2$, then there are two cases:
\begin{itemize}
\item $S$ is isogenous to product of curves of unmixed type $(C
\times F)/G$, and the classification of these surfaces is given by
the cases labelled with UnMix in Table 1;

\item $S$ is isogenous to a product of curves of mixed type $(C
\times C)/G$, there is only one such case and it is labelled with Mix
in Table 1;

\end{itemize}
\end{enumerate}
\begin{flushleft}
\begin{tabular}{|c|c|c|c|c|c|c|c|c|}
  \hline
 Type & $K^2_S$ & $g(F)$ & $g(C)$ & $G$ & \verb|IdSmallGroup| & \textbf{m} & $dim$ & $n$ \\
  \hline
 GH & $8$ & $2$ & $3$ &  $\mathbb{Z}/2\ZZ$ &                    G(2,1) & $(2^6)$ & $6$ & $1$ \\ \hline
 GH & $8$ &   $2$ & $4$ &  $\mathbb{Z}/3\ZZ$ &                     G(3,1) & $(3^4)$ & $4$ & $1$ \\ \hline
 GH & $8$ &   $2$ & $5$ &  $\mathbb{Z}/2\ZZ \times \mathbb{Z}/2\ZZ$ & G(4,2) & $(2^5)$ & $5$ & $2$  \\ \hline
 GH & $8$ &   $2$ & $5$ &  $\mathbb{Z}/4\ZZ$ &                     G(4,1) & $(2^2,4^2)$ & $4$ & $1$ \\ \hline
 GH & $8$ &   $2$ & $6$ &  $\mathbb{Z}/5\ZZ$ &                     G(5,1) & $(5^3)$ & $3$ & $1$ \\ \hline
 GH & $8$ &   $2$ & $7$ &  $\mathbb{Z}/6\ZZ$ &                     G(6,2)  & $(2^2,3^2)$ & $4$ & $1$\\ \hline
 GH & $8$ &  $2$ & $7$ &  $\mathbb{Z}/6\ZZ$ &                     G(6,2) & $(3,6^2)$ & $3$ & $1$\\ \hline
 GH & $8$ &   $2$ & $9$ &  $\mathbb{Z}/8\ZZ$ &                     G(8,1) & $(2,8^2)$ & $3$ & $1$\\ \hline
 GH & $8$ &  $2$ & $11$ & $\mathbb{Z}/{10}\ZZ$&                   G(10,2) & $(2,5,10)$ & $3$ & $1$ \\ \hline
 GH & $8$ &  $2$ & $13$ & $\mathbb{Z}/2\ZZ \times \mathbb{Z}/6\ZZ$ & G(12,5) & $(2,6^2)$ & $3$ & $2$\\ \hline
 GH & $8$ &  $2$ & $7$ &  $\mathfrak{S}_3$ &                              G(6,1) & $(2^2,3^2)$ & $4$ & $1$\\ \hline
 GH & $8$ &  $2$ & $9$ &  $Q_8$ &                              G(8,4) & $(4^3)$ & $3$ & $1$ \\ \hline
 GH & $8$ &  $2$ & $9$ &  $D_4$ &                              G(8,3) & $(2^3,4)$ & $4$ & $2$ \\ \hline
 GH & $8$ &  $2$ & $13$ & $D_6$ &                              G(12,4) & $(2^3,3)$ & $3$ & $2$ \\ \hline 
 GH & $8$ &  $2$ & $13$ & $D_{4,3,-1}$ &                       G(12,1) & $(3,4^2)$ & $3$ & $1$ \\ \hline
 GH & $8$ &  $2$ & $17$ & $D_{2,8,3}$ &                        G(16,8) & $(2,4,8)$ & $3$ & $1$ \\ \hline
 GH & $8$ &  $2$ & $25$ & $\mathbb{Z}/2\ZZ \ltimes ((\mathbb{Z}/2\ZZ)^2 \times \mathbb{Z}/3\ZZ)$ & G(24,8) & $(2,4,6)$ & $3$ & $2$ \\ \hline
 GH & $8$ &  $2$ & $25$ & $\SL(2,\mathbb{F}_3)$ &               G(24,3) & $(3^2,4)$ & $3$ & $1$ \\ \hline
 GH & $8$ &  $2$ & $49$ & $\GL(2,\mathbb{F}_3)$ &               G(48,29) & $(2,3,8)$ & $3$ & $1$ \\ \hline
UnMix & $8$ & $3$ & $3$ & $\mathbb{Z}/2\ZZ \times \mathbb{Z}/2\ZZ$
& G(4,2) & $(2^2)$, $(2^2)$ & $4$ &$1$ \\ \hline UnMix & $8$ & $3$
& $4$ & $\mathfrak{S}_3$                              & G(6,1) & $(3)$,
$(2^2)$ & $3$ &$1$  \\ \hline UnMix & $8$ &  $3$ & $5$ & $D_4$ &
G(8,3) & $(2)$, $(2^2)$ & $3$ & $1$ \\ \hline Mix & $8$ & $3$ &
$3$  & $\mathbb{Z}/4\ZZ$ &                      G(4,1) & - & $3$ &
$1$
\\ \hline
\end{tabular}
\end{flushleft}
\begin{center}
Table 1.
\end{center}

In Table 1 \verb|IdSmallGroup| denotes the label of the group
$G$ in the GAP4 database of small groups, $\mathbf{m}$ is the
branching data, i.e., the number of points (exponent) of given branching order (base) of the covering $F \rightarrow \mathbb{P}^1$ (e.g. $2^6$ means that there are six points with ramification order two). Moreover, each item provides a union of connected
components of the moduli space of surfaces of general type, their
dimension is listed in the column $dim$, and $n$ is the number of
connected components.

\end{theo}

The proof of the theorem involves techniques from both geometry
and combinatorial group theory developed in \cite{BCG08} and
\cite{CP09}.

An idea of the proof of the theorem is the following. There are
two cases. In the first case $S$ is not of Albanese maximal dimension; 
it is of generalized hyperelliptic type, hence one
classifies all possible finite groups $G$ which induce a
$G-$covering $\pi_F\colon F \rightarrow \mathbb{P}^1$ with
$g(F)=2$ (see \cite{bolza}). Then we need to check whether such groups $G$ induce
an unramified $G-$covering $\pi_C\colon C \rightarrow B \cong
C/G$, where the genus of $B$ is $2$ and the genus of $C$ is
determined by the Riemann-Hurwitz formula. We notice that the
action of $G$ on the product $C \times F$ is always free, because
the action on $C$ is free.

In the second case, $S$ is of maximal Albanese dimension and
it is slightly more difficult. In this case, both projections $C
\rightarrow C/G$ and $C \rightarrow F/G$ are ramified covers of
elliptic curves. By  \cite[Corollary 2.4]{zucc} the genera of $C$
and $F$ can be at most five. Therefore, we can proceed as in the
previous case by analyzing all possible combinations of genera. This
time we also take into account the fact that, since we do not have
an \'{e}tale cover, we have to check whether the action of $G$ on
the product of the two curves is free or not.

\begin{rem} By the above theorem, we know that
$\mathcal{M}_{8,2,2}$ consists of at least $28$ connected
components of dimensions between $3$ and $6$.
\end{rem}

An interesting and still unanswered question is the following one

\begin{quest}
\textit{Is there a surface with $p_g=q=2$ and $K^2_S=8$ which is
not isogenous to a product?}
\end{quest}

It is worth noticing that up to now the only known examples of
surfaces of general type with $\chi(S)=1$ and $K^2_S=8$ are all
isogenous to a higher product of curves. These surfaces have
all been classified, in particular for $p_g=q=1$ see \cite{CP09},
and for $p_g=q=0$ see \cite{BCG08}.

\medskip

The study of surfaces of general type with $p_g=q=2$ and
$\textrm{Albdim}(S)=2$ started with a generalization of the
definition of surface isogenous to a product. A surface $S$ is
said to be a \emph{product-quotient} surface if  $S$ is the
minimal model of the minimal desingularization of the quotient
surface $T:=(C_1 \times C_2)/G$ where $C_1$, $C_2$ are curves of
genus at least two, and $G$ is a finite group acting diagonally
but not necessarily freely on the product.
\begin{theo}[\cite{Pe11} Theorem 1.1]\label{theo.prod_quot}
Let $S$ be a product-quotient surface of general type with
$p_g=q=2$. In particular $S$ is the minimal desingularization of
$T:=(C_1 \times C_2)/G$, and these surfaces are classified in
Table 2.

\begin{center}
\begin{tabular}{|c|c|c|c|c|c|c|c|c|c|}
  \hline
 $K^2_S$ & $g(C_1)$ & $g(C_2)$ & $G$ & \verb|IdSmallGroup| & \textbf{m} & Type & Num.~Sing.  & dim & n
\\ \hline
$4$ &  $2$ & $2$ &  $\mathbb{Z}/2\ZZ$  & G(2,1)  & $(2^2)$ $(2^2)$&
$\frac{1}{2}(1,1)$ & $4$ & $4$ & $1$
\\ \hline $4$ &  $3$ & $3$ &  $D_4$ & G(8,3)  & $(2)$ $(2)$ & $\frac{1}{2}(1,1)$ & $4$ & $2$ & $1$
\\ \hline $4$ & $3$ & $3$ &  $Q_8$ & G(8,4)  & $(2)$ $(2)$ & $\frac{1}{2}(1,1)$ & $4$ & $2$ & $1$
\\ \hline $5$ &  $3$ & $3$ &  $\mathfrak{S}_3$ & G(6,1) &  $(3)$ $(3)$ & $\frac{1}{3}(1,1)+\frac{1}{3}(1,2)$ & $2$ & $2$ & $1$
\\ \hline $6$ &  $4$ & $4$ &  $\mathfrak{A}_4$ & G(12,3) & $(2)$ $(2)$ & $\frac{1}{2}(1,1)$ & $2$ & $2$ & $1$
\\
\hline
\end{tabular}
Table 2.
\end{center}
In Table 2 each item provides a union of irreducible subvarieties
of the moduli space of surfaces of general type, their dimension
is listed in the column $dim$ and $n$ is the number of
subvarieties. Moreover the columns of Table 2 labelled with $Type$
and $Num.~Sing.$ indicate the types and the number of
singularities of $T$.
\end{theo}


We observe that a complete classification of product-quotient
surfaces of general type with $\chi(S)=1$ and $p_g \leq 1$ is
still missing. A partial classification for the ones with
$p_g=q=1$ is given in \cite{MP10} and for the ones with $p_g=q=0$
in \cite{BP11}.
%

\section{Surfaces with $p_g=q=2$ and $K^2_S=6$}

%
\subsection*{Albanese map of degree 4}

We briefly recall the construction of the surface with $p_g=q=2$
and $K^2_S=6$ given in Theorem \ref{theo.prod_quot}. Let us denote
by $\mathfrak{A}_4$ the alternating group on four symbols and by
$V_4$ its Klein subgroup, namely
\begin{equation*}
V_4=\langle id, \, (12)(34), \, (13)(24), \, (14)(23) \rangle
\cong (\mathbb{Z}/2 \mathbb{Z})^2.
\end{equation*}
$V_4$ is normal in $\mathfrak{A}_4$ and  the quotient
$H:=\mathfrak{A}_4 / V_4$ is a cyclic group of order $3$. By
Riemann's existence theorem, it is possible to construct two smooth
curves $C_1$, $C_2$ of genus $4$ endowed with an action of
$\mathfrak{A}_4$ such that the only non-trivial stabilizers are
the elements of $V_4$. Then
\begin{itemize}
\item $E_i':=C_i/\mathfrak{A}_4$ is an elliptic curve; \item the
$\mathfrak{A}_4$-cover $f_i \colon C_i \to E_i'$ is branched at
exactly one point of $E_i'$, with branching order $2$.
\end{itemize}
It follows that the surface
\begin{equation*}
\widehat{X}:=(C_1 \times C_2)/\mathfrak{A}_4,
\end{equation*}
where $\mathfrak{A}_4$ acts diagonally, has  two rational double
points of type $\frac{1}{2}(1, \, 1)$ and no other singularities.
It is straightforward to check that the desingularization
$\widehat{S}$ of $\widehat{X}$ is a minimal surface of general
type with $p_g=q=2$, $K_{\widehat{S}}^2=6$ and that $\widehat{X}$ is the
canonical model of $\widehat{S}$.

The $\mathfrak{A}_4$-cover $f_i \colon C_i \to E_i'$ factors
through the bidouble cover $g_i \colon C_i \to E_i$, where
$E_i:=C_i/V_4$. Note that $E_i$ is again an elliptic curve, so
there is an isogeny $E_i \to E_i'$, which is a triple Galois cover
with Galois group $H$. Consequently, we have an isogeny
\begin{equation*}
p \colon E_1 \times E_2 \to \widehat{A}:=(E_1 \times E_2)/H,
\end{equation*}
where the group $H$ acts diagonally, and a commutative diagram
\begin{equation} \label{dia.example.pq}
\begin{xy}
\xymatrix{
C_1 \times C_2  \ar[d]_{g_1 \times g_2} \ar[rr]^{\tilde{p}} & & (C_1 \times C_2 )/ \mathfrak{A}_4 =\widehat{X} \ar[d]^{\widehat{\alpha}}  \\
E_1 \times E_2   \ar[rr]^{p} & & (E_1 \times E_2)/H = \widehat{A}.   \\
  }
\end{xy}
\end{equation}
The morphism $\hat{\alpha} \colon \widehat{X} \to \widehat{A}$ is
the Albanese map of $\widehat{X}$; it is a finite, non-Galois
quadruple cover.

This last observation led us to study surfaces $\widehat{S}$ with
$p_g=q=2$ and $K^2_{\widehat{S}}=6$ whose Albanese map
$\hat{\alpha} \colon \widehat{S} \to
\widehat{A}:=\textrm{Alb}(\widehat{S})$ is a quadruple cover of an
abelian surface $\widehat{A}$, in order to understand and describe
the deformations of $\widehat{S}$.

The main result in the theory of quadruple cover for algebraic
varieties is the following.

\begin{theo}\emph{\cite[Theorem 1.2]{HM99}}\label{theo.mir}
Let $Y$ be a smooth algebraic variety. Any quadruple cover
$f\colon X \longrightarrow Y$ is determined by a locally free
$\mathcal{O}_Y$-module $\mathcal{E}^{\vee}$ of rank $3$ and a
totally decomposable section $\eta \in H^0(Y, \, \bigwedge^2 S^2
\mathcal{E}^{\vee} \otimes \bigwedge^3 \mathcal{E})$.
\end{theo}

The vector bundle $\EE^{\vee}$ is called the \emph{Tschirnhausen
bundle} of the cover. 
\\ We have $f_{*}\mathscr{O}_X= \mathscr{O}_Y
\oplus \EE^{\vee}$.

In the situation above, we can prove that the Tschirnhausen bundle
 $\mathcal{E}^{\vee}$ associated with the Albanese cover
 is of the form
$\Phi^{\mathcal{P}}(\mathcal{L}^{-1})^{\vee}$, where $\mathcal{L}$ is a
polarization on $A$ (the dual abelian variety of $\wA$) and
$\Phi^{\mathcal{P}}$ denotes the Fourier-Mukai transform with kernel the normalized Poincar\'e bundle $\mathcal{P}$ i.e.,
\begin{equation*}
\mathcal{E}=\Phi^{\mathcal{P}}(\mathcal{L}^{-1}):={R}^i\pi_{\wA \,
*}(\mathcal{P} \otimes \pi^*_{A}\mathcal{L}^{-1}).
\end{equation*}
 More
precisely, the bundle $\mathscr{E}^{\vee}$ has rank $3$ and
$\mathscr{L}$ is a polarization of type $(1, \, 3)$.

We are able to prove the following theorem.

\begin{theo}\emph{\cite[Theorem 2.1]{PP12}}
There exists a $4$-dimensional family $\mathcal{M}_{\Phi}$ of
minimal surfaces of general type with $p_g=q=2$ and $K^2_{\widehat{S}}=6$ such
that, for the general element $\widehat{S} \in
\mathcal{M}_{\Phi}$, the canonical class $K_{\wS}$ is ample and
the Albanese map $\hat{\alpha}\colon \widehat{S} \longrightarrow
\widehat{A}$ is a finite cover of degree $4$.

The Tschirnhausen bundle $\EE^{\vee}$ associated with
$\hat{\alpha}$ is isomorphic to
$\Phi^{\mathscr{P}}(\mathscr{L}^{-1})$, where $\mathscr{L}$ is a
polarization of type $(1, \, 3)$ on $A$.

The family $\mathcal{M}_{\Phi}$ provides an irreducible component
of the moduli space $\mathcal{M}_{6, \, 2, \,2}$ of canonical
models of minimal surfaces of general type with $p_g=q=2$ and
$K^2=6$. Such a component is generically smooth and contains the
$2$-dimensional family of product-quotient surfaces given in
Theorem \ref{theo.prod_quot}.
\end{theo}

This theorem is obtained by extending the construction given in
\cite{CH06} to the much more complicated case of quadruple covers.
More precisely, in order to build a quadruple cover $\hat{\alpha}
\colon \widehat{S} \lr \widehat{A}$ with Tschirnhausen bundle
$\EE^{\vee}$, we first build a quadruple cover $\alpha \colon S
\lr A$ with Tschirnhausen bundle $\phi_{\LL^{-1}}^*\EE^{\vee}$
(here $\phi_{\LL} \colon A \to \widehat{A}$ denotes the group
homomorphism sending $x \in A$ to $t_x^* \mathscr{L}^{-1} \otimes
\mathscr{L} \in \widehat{A}$) and then, by using the Schr\"odinger
representation of the finite Heisenberg group $\mathscr{H}_3$ on
$H^0(A, \, \LL)$, we identify the covers of this type which
descend to a quadruple cover $\hat{\alpha} \colon \widehat{S} \lr
\widehat{A}$. For the general surface $\wS$, the branch divisor
$\widehat{B} \subset \wA$ of $\alpha \colon \wS \lr \wA$ is a
curve in the linear system $|\widehat{\mathscr{L}}^{\otimes 2}|$,
where $\widehat{\mathscr{L}}$ is a polarization of type $(1, \,3)$
on $\wA$. $\widehat{B}$ has six ordinary cusps and no other
singularities; such a curve is $\mathscr{H}_3$-equivariant and can
be associated with the dual of a member of the Hesse pencil of
plane cubics in $\mathbb{P}^2$. Once we have established these
facts with the help of a \verb|MAGMA| script, we can show that
$\widehat{S}$ is smooth in the general case.


\subsection*{Albanese map of degree 2}

In \cite{PP11} we give other families of surfaces with
$p_g=q=2$, $K^2_S=6$, and with Albanese map of degree equal to
two.

Let us introduce some notation. Let $(A, \,\mathcal{L})$ be a $(1,
\,2)$-polarized abelian surface and let us denote by $\phi_2
\colon A[2] \to \widehat{A}[2]$ the restriction of the canonical
homomorphism $\phi_{\mathcal{L}} \colon A \to \widehat{A}$ to the
subgroup of $2$-division points. Then $\textrm{im}\, \phi_2$
consists of four line bundles $\{\mathcal{O}_A, \, \mathcal{Q}_1,
\, \mathcal{Q}_2, \, \mathcal{Q}_3\}$. Let us denote by
$\textrm{im}\, \phi_2^{\,\times}$ the set $\{\mathcal{Q}_1, \,
\mathcal{Q}_2, \, \mathcal{Q}_3\}$.

\begin{theo}\emph{\cite[Theorem A]{PP11}} \label{theo.A}
Given an  abelian surface $A$ with a symmetric polarization $\mL$
of type $(1, \,2)$, not of product type,  for any $\mQ \in
\emph{im}\, \phi_2$ there exists a curve $D \in |\mL^2 \otimes
\mQ|$ whose unique singularity is an ordinary quadruple point at
the origin $o \in A$. Let $\mQ^{1/2}$ be a square root of $\mQ$,
and if $\mQ=\oo_A$ assume moreover $\mQ^{1/2} \neq \oo_A$. Then
the minimal desingularization $S$ of the double cover of $A$
branched over $D$ and defined by $\mL \otimes \mQ^{1/2}$ is a
minimal surface of general type with $p_g=q=2$, $K_S^2=6$ and
Albanese map of degree $2$.

Conversely, every minimal surface of general type with $p_g=q=2$,
$K_S^2=6$ and Albanese map of degree $2$ can be constructed in
this way.
\end{theo}

When $\mQ=\mQ^{1/2}=\oo_A$, we obtain a minimal surface
with $p_g=q=3$, which correspond to the second case of Theorem
\ref{theo_pgq3}.

We use the following terminology:
\begin{itemize}
\item if $\mQ =\oo_A$ we say that $S$ is a \emph{surface of type}
$I$. Furthermore, if $\mathcal{Q}^{1/2} \notin \textrm{im} \,
\phi_2$ we say that $S$ is of type $Ia$, whereas if
$\mathcal{Q}^{1/2} \in \textrm{im} \, \phi_2^{\,\times}$  we say
that $S$ is of type $Ib;$ \item if $\mQ \in
\textrm{im}\,\phi_2^{\, \times}$ we say that $S$ is a
\emph{surface of type} $II$.
\end{itemize}

Recall that the degree of the Albanese map in this case is a
topological invariant since $q(S)=2$. Therefore, let us consider the moduli space
$\mathcal{M}$ of minimal surfaces of general type with $p_g=q=2$,
$K_S^2=6$ and Albanese map of degree $2$. Let $\mathcal{M}_{Ia}$,
$\mathcal{M}_{Ib}$, $\mathcal{M}_{II}$ be the algebraic subsets
whose points parametrize isomorphism classes of surfaces of type
$Ia$, $Ib$, $II$, respectively. Therefore $\mathcal{M}$ can be
written as the disjoint union
\begin{equation*}
\mathcal{M}=\mathcal{M}_{Ia} \sqcup \mathcal{M}_{Ib} \sqcup
\mathcal{M}_{II}.
\end{equation*}

Our result on the moduli space is

\begin{theo}\emph{\cite[Theorem B]{PP11}} \label{theo.B} The following holds:
\begin{itemize}
\item[$\boldsymbol{(i)}$] $\mathcal{M}_{Ia}$, $\mathcal{M}_{Ib}$,
$\mathcal{M}_{II}$ are the connected components of $\mathcal{M};$
\item[$\boldsymbol{(ii)}$] these are also \emph{irreducible}
components of the moduli space of minimal surfaces of general
type$;$ \item[$\boldsymbol{(iii)}$] $\mathcal{M}_{Ia}, \,
\mathcal{M}_{Ib}, \, \mathcal{M}_{II}$ are generically smooth, of
dimension $4, \, 4, \, 3$, respectively$;$
\item[$\boldsymbol{(iv)}$] the general surface in
$\mathcal{M}_{Ia}$ and $\mathcal{M}_{Ib}$ has ample canonical
class$;$ all surfaces in $\mathcal{M}_{II}$ have ample canonical
class.
\end{itemize}
\end{theo}

Some interesting questions remain open also for surfaces with $p_g=q=2$ and $K^2_S=6$ namely

\begin{quest}
\begin{itemize}

\item Is $\mathcal{M}_{\Phi}$ a \emph{connected} component of
$\mathcal{M}_{6, \, 2,\, 2}?$

\item What are the possible degrees for the Albanese map of a
minimal surface with $p_g=q=2$ and $K^2=6$?

\item What are the irreducible/connected components of
$\mathcal{M}_{6, \, 2,\, 2}$?

\end{itemize}
\end{quest}
%

\section{Surfaces with $p_g=q=2$ and $K^2_S=5$}

A similar calculation as the one in the previous section
shows that any element of the family of surfaces with $p_g=q=2$ and
$K^2_S=5$ given in Theorem \ref{theo.prod_quot} has Albanese map
of degree $3$. Moreover, for such element $S$ the image $\alpha(S)$ is an abelian surface with a
polarization of type $(1,2)$.
In order to understand the connected components of the moduli space that contains this family, we have to study triple covers in detail. 


Let us introduce some terminology. Let $S$ be a minimal surface of
general type with $p_g=q=2$ and $K_S^2=5$, such that its Albanese
map $\alpha \colon S \to \textrm{Alb}(S)$
 is a generically finite morphism of degree $3$. If we consider the Stein
factorization of $\alpha$, i.e.,
\begin{equation*}
S \stackrel{p}{\lr} \widehat{X} \stackrel{\hat{f}}{\lr}
\textrm{Alb}(S),
\end{equation*}
the map $\hat{f} \colon \widehat{X} \to \textrm{Alb}(S)$ is a
flat triple cover, which can be studied by applying the techniques
developed in \cite{M85}. In particular, $\hat{f}$ is determined by
a rank $2$ vector bundle $\EE^{\vee}$ on $\textrm{Alb}(S)$, called
the \emph{Tschirnhausen bundle} of the cover, and by a global
section $\eta \in H^0(\textrm{Alb}(S), \, S^3\EE^{\vee} \otimes
\bigwedge^2 \EE)$. In the examples of \cite{CH06} and \cite{Pe11},
the surface $\widehat{X}$ is singular; nevertheless, in both cases
the numerical invariants of $\EE$ are the same predicted by the
formulae of \cite{M85}, as if $\widehat{X}$ were smooth. This
leads us to introduce the definition of \emph{negligible
singularity} for a triple cover, which is similar to the
corresponding well-known definition for double covers. Then,
inspired by the construction in \cite{CH06}, we say that $S$ is a
\emph{Chen-Hacon surface} if there exists a polarization
$\mathcal{L}$ of type $(1,\, 2)$ on
$\textrm{Pic}^0(S)=\widehat{\textrm{Alb}(S)}$ such that
$\EE^{\vee}$ is the Fourier-Mukai transform of the line bundle
$\mathcal{L}^{-1}$.

Our first result is the following characterization of
 Chen-Hacon surfaces.

\begin{theo}\emph{\cite[Theorem A]{PP10}}
Let $S$ be a minimal surface of general type with $p_g=q=2$ and
$K_S^2=5$ such that the Albanese map $\alpha \colon S \to
\textrm{Alb}(S)$ is a generically finite morphism of degree $3$.
Let
\begin{equation*}
S \stackrel{p}{\lr} \widehat{X} \stackrel{\hat{f}}{\lr}
\textrm{Alb}(S)
\end{equation*}
be the Stein factorization of $\alpha$. Then $S$ is a Chen-Hacon
surface if and only if $\widehat{X}$ has only negligible
singularities.
\end{theo}

Moreover, we can completely describe all the possibilities for the
singular locus of $\widehat{X}$. It follows that $\widehat{X}$ is never
smooth, and it always contains a cyclic quotient singularity of
type $\frac{1}{3}(1, \, 1)$. Therefore $S$ always contains a
$(-3)$-curve, which turns out to be the fixed part of the
canonical system $|K_S|$.

Now let $\mathcal{M}_{5,2,2}$ be the moduli space of surfaces with
$p_g=q=2$ and $K^2_S=5$, and let $\mathcal{M}^{CH} \subset \mathcal{M}_{5,2,2}$ be the
subset of the points parametrizing (isomorphism classes of)
 Chen-Hacon surfaces. Our second result is the
following.

\begin{theo}\emph{\cite[Theorem B]{PP10}}\label{theo_B_k5}
$\mathcal{M}^{CH}$ is an irreducible, connected, generically
smooth component of $\mathcal{M}_{5,2,2}$ of dimension $4$.
\end{theo}


It is worth noticing that Theorem \ref{theo_B_k5} shows that every small deformation
of a Chen-Hacon surface is still a Chen-Hacon surface; in
particular, no small deformation of $S$ makes the $(-3)$-curve
disappear. Moreover, since $\mathcal{M}^{CH}$ is generically
smooth, the same is true for the first-order deformations. By
contrast, Burns and Wahl proved in \cite{BW74} that first-order
deformations always smooth all the $(-2)$-curves. Furthermore, Catanese
used this fact in \cite{Ca89} to produce examples of
surfaces of general type with everywhere non-reduced moduli
spaces. Theorem \ref{theo_B_k5} shows rather strikingly that the results
of Burns-Wahl and Catanese cannot be extended to the case of
$(-3)$-curves and, as far as we know, provides the first explicit
example of this situation.

\begin{quest}\label{qst.k5}
There are mainly three questions open for surfaces with $p_g=q=2$ and $K^2_S=5$:
\begin{itemize}
\item Are there surfaces with these invariants whose Albanese map
has degree different from $3$?
\item Are there surfaces with these
invariants whose Albanese map has degree $3$, but which are not
Chen--Hacon surfaces?
\item What are the irreducible/connected
components of $\mathcal{M}_{5, \, 2,\, 2}$?
\end{itemize}
\end{quest}

%

\section{Surfaces with $p_g=q=2$ and $K^2_S=4,7$ and $9$}

Surfaces with $p_g=q=2$ and $K^2_S=4$ were studied by
Ciliberto and Mendes Lopes in \cite{CML02}. The authors
classified all surfaces with $p_g=q=2$ and non-birational
bicanonical map (not presenting the standard case). Their result
is the following.
\begin{theo}
If $S$ is a minimal surface of general type with $p_g=q=2$ and
non-birational bicanonical map not presenting the standard case,
then $S$ is a double cover of a principally polarized abelian
surface $(A,\Theta)$, with $\Theta$ irreducible. The double cover
$S \rightarrow A$ is branched along a divisor $B \in |\Theta|$,
having at most double points. In particular $K^2_S =4$.
\end{theo}

Surfaces with $p_g=q=2$ and $K^2_S=4$ were also considered by Manetti in his work on Severi's conjecture
(\cite{M}). In particular, Manetti proved that if $p_g=q=2$, $K_S$ is \emph{ample}, and $K^2_S=4$, then $S$ is a double cover of its Albanese variety,
which turns out to be  a principally polarized abelian surface
$(A,\Theta)$. In our investigation on product-quotient surfaces in
\cite{Pe11}, we constructed three families of surfaces with
$p_g=q=2$ and $K^2_S=4$: see Theorem \ref{theo.prod_quot}. In all
these cases, the canonical bundle is not ample; nevertheless, the
Albanese map has degree two and the Albanese variety is a
principally polarized abelian surface. 
\\ Eventualy in \cite{CMLP13}, it is proven that all surfaces with $p_g=q=2$ and $K^2_S=4$ are double covers of a principally polarized abelian surface. This also proves that there is only one connected component of the moduli space $\mathcal{M}_{4,2,2}$ of dimension $4$.

The cases with $K^2_S=7$ and $9$ are more mysterious. There are no
examples known up to now.

As we have seen, the degree of the Albanese map plays a very
important role in the classification of surfaces with $p_g=q=2$,
so it would be nice to have an upper bound on this number. We
notice that the degree of the Albanese map of all the families given in Table 0 is always less or equal
to $K^2_S -2$. Unfortunately we cannot prove that this is the
actual upper bound.
%

%

\bigskip
\bigskip



\begin{thebibliography} {9}
%
\bibitem[BCP06]{BCP06}
I. Bauer, F. Catanese, R. Pignatelli, \textit{Complex surfaces of
general type: some recent progress}. Global aspects of complex
geometry. Springer, Berlin (2006), 1--58.
%
\bibitem[BCG08]{BCG08} 
I. Bauer, F. Catanese, F. Grunewald, \textit{The classification of
surfaces with $p_g=q=0$ isogenous to a product}. Pure Appl. Math.
Q., \textbf{4}, no. 2, part1, (2008), 547--586.
%
\bibitem[BP11]{BP11}
I. Bauer, R. Pignatelli, \textit{The classification of minimal
product-quotient surfaces with $p_g=0$}. Mathematics of
Computation \textbf{81}, 280 (2012), 2389--2418.
%
\bibitem[B78]{Be} 
A. Beauville, \textit{Surfaces Alg\'{e}briques Complex}.
Ast\'{e}risque \textbf{54}, Paris (1978).
%
\bibitem[B82]{Be82} 
A. Beauville, \textit{L'inegalite $p_g \geq 2q-4$ pour les surface
de type general}. Bull. Soc. Math. France \textbf{110}, (1982),
344--346.
%
%
\bibitem[Bol88]{bolza}
O. Bolza, \textit{On binary sextic with linear transformations into
themselves}. Amer. J. Math. \textbf{10}, (1888), 47--70.
%
%
\bibitem[BW74]{BW74}
D. M. Burns, J. M. Wahl, \textit{Local contributions to global
deformations of surfaces}. Invent. Math. \textbf{26} (1974), 67-88.
%
\bibitem[CP09]{CP09}
G. Carnovale, F. Polizzi, \textit{The classification of surfaces
with $p_g=q=1$ isogenous to a product of curves}. Adv. Geom.,
Vol. \textbf{9}, no. 2, (2009), 233--256.
%
\bibitem[Ca89]{Ca89}
F. Catanese, \textit{Everywhere non reduced moduli spaces}.
Invent. Math. \textbf{98} (1989), 293-310.
%
\bibitem[Ca91]{Ca91}
F. Catanese, \textit{Moduli and classification of irregular
Kaehler manifolds (and algebraic varieties) with Albanese general
type fibrations}. Invent. Math. \textbf{104} (1991), 263-289.
%
\bibitem[Ca00]{cat00} 
F. Catanese, \textit{Fibred surfaces, varieties isogenous to a
product and related moduli spaces}. Amer. J. Math. \textbf{122},
(2000), 1--44.
%
\bibitem[Ca11]{Ca11}
F. Catanese, \textit{A superficial working guide to deformations
and moduli}, e-print arXiv:1106.1368, to appear in in the Handbook
of Moduli, a volume in honour of David Mumford, to be published by
International Press.
%
\bibitem[CCML98]{CCML}
F. Catanese, C. Ciliberto, M. Mendes Lopes, \textit{On the
classification of irregular surfaces of general type with
nonbirational bicanonical map}. Trans. Amer. Math. Soc.
\textbf{350}, no. 1, (1998), 275--308.
%
\bibitem[CH06]{CH06}
J. Chen, C. Hacon, \textit{A Surface of general type with
$p_g=q=2$ and $K^2=5$}. Pacific. J. Math. \textbf{223} no. 2,
(2006), 219--228.
%
\bibitem[CML02]{CML02}
C. Ciliberto, M. Mendes Lopes, \textit{On surfaces with $p_g=q=2$
and non-birational bicanonical map}. Adv. Geom. \textbf{2} n. 3,
(2002), 281--300.
%
\bibitem[CMLP13]{CMLP13}
C. Ciliberto, M. Mendes Lopes, R. Pardini, \textit{The classification of minimal irregular surfaces of general type with $K^2= 2p_g$}, preprint arXiv:1307.6228
%
\bibitem[Gie77]{Gie77}
D. Gieseker, \textit{Global moduli for surfaces of general type}.
Invent. Math \textbf{43} no. 3, (1977), 233--282.
%
\bibitem[HP02]{HP02}
C. Hacon, R. Pardini, \textit{Surfaces with $p_g=q=3$}. Trans.
Amer. Math. Soc. \textbf{354}, (2002), 2631--2638.
%
\bibitem[HM99]{HM99}
D.W. Hahn, R. Miranda, \textit{Quadruple covers of algebraic
varieties}. J. Algebraic Geom. \textbf{8}, 1999.
%
\bibitem[M03]{M}
M. Manetti, \textit{Surfaces of Albanese General type and the
Severi conjecture}. Math. Nachr., Vol. \textbf{261/262}, (2003),
105--122.
%
\bibitem[M85]{M85}
R. Miranda, \textit{Triple covers in algebraic geometry}. Amer. J.
of Math. \textbf{107} (1985), 1123-1158.
%
\bibitem[MP10]{MP10}
E. Mistretta, F. Polizzi, \textit{Standard isotrivial fibrations
with $p_g=q=1$ II}. AG/0805.1424.
%
\bibitem[Pi02]{Pi02}
G.P. Pirola, \textit{Surfaces with $p_g=q=3$}. Manuscripta Math.
\textbf{108} no. 2, (2002), 163--170.
%
\bibitem[Pe11]{Pe11}
M. Penegini, \textit{The classification of isotrivially fibred
surfaces with $p_g=q=2$}. Collect. Math. \textbf{62}, No.
\textbf{3}, (2011), 239--274.
%
\bibitem[PP10]{PP10}
M. Penegini, F. Polizzi, \textit{On surfaces with $p_g=q=2$,
$K^2=5$ and Albanese map of degree $3$}. e-print arXiv:1011.4388
(2010), to appear in Osaka J. Math.
%
\bibitem[PP11]{PP11}
M. Penegini, F. Polizzi, \textit{On surfaces with $p_g=q=2$,
$K^2=6$ and Albanese map of degree $2$}. e-print arXiv:1105.4983
(2011), to appear in Canad. J. Math.
%
\bibitem[PP12]{PP12}
M. Penegini, F. Polizzi, \textit{A new family of surfaces with
$p_g=q=2$ and $K^2=6$ whose Albanese map has degree $4$}. e-print
arXiv:1207.3526.
%
\bibitem[Z03]{zucc}
F. Zucconi, \textit{Surfaces with $p_g=q=2$ and an irrational
pencil}. Canad. J. Math. Vol. \textbf{55}, (2003), 649--672.


\end{thebibliography}
\end{document}